\documentclass[12pt]{amsart}
\usepackage{graphicx}
\usepackage{amssymb}
\usepackage{amsfonts}
\usepackage{amsmath}
\usepackage{array}
\usepackage{pstricks}
\usepackage{tableaux}
\usepackage{bm}
\usepackage[utf8]{inputenc}

\usepackage{pspicture}
\usepackage[matrix,ps,xdvi]{xypic} 

\newdimen\Squaresize \Squaresize=14pt
\newdimen\Thickness \Thickness=0.5pt

\def\Square#1{\hbox{\vrule width \Thickness
   \vbox to \Squaresize{\hrule height \Thickness\vss
      \hbox to \Squaresize{\hss#1\hss}
   \vss\hrule height\Thickness}
\unskip\vrule width \Thickness}
\kern-\Thickness}

\def\Vsquare#1{\vbox{\Square{$#1$}}\kern-\Thickness}

\def\young#1{
\vbox{\smallskip\offinterlineskip
\halign{&\Vsquare{##}\cr #1}}}

\def\boxit#1#2{\setbox1=\hbox{\kern#1{#2}\kern#1}%
\dimen1=\ht1 \advance\dimen1 by #1 \dimen2=\dp1 \advance\dimen2 by #1
\setbox1=\hbox{\vrule height\dimen1 depth\dimen2\box1\vrule}%
\setbox1=\vbox{\hrule\box1\hrule}%
\advance\dimen1 by .4pt \ht1=\dimen1
\advance\dimen2 by .4pt \dp1=\dimen2 \box1\relax}
\def\bo#1{\boxit{1pt}{$#1$}}

\newdimen\vcadre\vcadre=0.1cm 
\newdimen\hcadre\hcadre=0.1cm 

\def\arx#1[#2]{\ifcase#1 \relax \or%
  \ar @{-}[#2]  \or%
  \ar @2{-}[#2] \or%
  \ar @{--}[#2] \or%
  \ar @2{.}[#2] \or%
  \ar @{~}[#2]  \fi}

\newtheorem{example}{Example}[section]
\newtheorem{note}[example]{Note}
\newtheorem{theorem}[example]{Theorem}

\newtheorem{corollary}[example]{Corollary}

\newtheorem{proposition}[example]{Proposition}

\newtheorem{lemma}[example]{Lemma}

\def\Proof{\noindent \it Proof -- \rm}
\def\qed{\hspace{3.5mm} \hfill \vbox{\hrule height 3pt depth 2 pt width 2mm}
\bigskip}

\def\alph{{\rm alph}}

\def\qbin#1#2{\left[ \begin{matrix} #1 \\ #2 \end{matrix} \right]_q}

\def\btr{\blacktriangleright}

\def\FQSym{{\bf FQSym}}

\def\WQSym{{\bf WQSym}}
\def\WSym{{\bf WSym}}
\def\Sym{{\bf Sym}}
\def\NCSF{{\bf Sym}}
\def\QSym{{\it QSym}}

\def\PBT{{\bf PBT}}

\def\ssh{\Cup}

\def\Std{{\rm Std}}

\def\<{\langle}
\def\>{\rangle}

\def\K{\operatorname{\mathbb K}}

\def\F{{\bf F}}

\def\G{{\bf G}}
\def\M{{\bf M}}
\def\P{{\bf P}}

\def\SG{{\mathfrak S}}

\def\shuff#1#2{\mathbin{
\hbox{\vbox{ \hbox{\vrule \hskip#2 \vrule height#1 width 0pt
}%
\hrule}%
\vbox{ \hbox{\vrule \hskip#2 \vrule height#1 width 0pt
\vrule }%
\hrule}%
}}}

\def\shuf{{\mathchoice{\shuff{7pt}{3.5pt}}%
{\shuff{6pt}{3pt}}%
{\shuff{4pt}{2pt}}%
{\shuff{3pt}{1.5pt}}}}%
\def\shuffle{\,\shuf\,}

\def\B{{\bm B}}

\def\maj{{\rm maj}}
\def\imaj{{\rm imaj}}
\def\inv{{\rm inv}}

\def\std{{\rm std}}

\title[]%
{Noncommutative Bell polynomials and the dual immaculate basis}
\author[J.-C.~Novelli, J.-Y.~Thibon, F. Toumazet]%
{Jean-Christophe Novelli, Jean-Yves Thibon and Fr\'ed\'eric Toumazet}

\address[]{[Novelli, Thibon, Toumazet] Laboratoire d'informatique Gaspard-Monge\\
Universit\'e Paris-Est Marne-la-Vall\'ee \\
5, Boulevard Descartes \\ Champs-sur-Marne \\
77454 Marne-la-Vall\'ee cedex 2 \\
France}
\email[Jean-Christophe Novelli]{novelli@u-pem.fr}
\email[Jean-Yves Thibon]{jyt@u-pem.fr} 
\email[ Fr\'ed\'eric Toumazet]{frederic.toumazet@u-pem.fr}
\date{\today}


\date{\today}

\begin{document}

\begin{abstract}
We define a new family of noncommutative Bell polynomials in the algebra of
free quasi-symmetric functions and relate it to the dual immaculate basis of
quasi-symmetric functions. We obtain noncommutative versions of Grinberg's
results [Canad. J. Math. 69 (2017), 21--53], and interpret these in terms of
the tridendriform structure of $\WQSym$. We then present a variant of Rey's
self-dual Hopf algebra of set partitions [FPSAC'07, Tianjin] adapted to our
noncommutative Bell polynomials and give a complete description of the Bell
equivalence classes as linear extensions of explicit posets.
\end{abstract}

\maketitle

\section{Introduction}

Many classical combinatorial numbers or polynomials are the dimensions (or
graded dimensions) of certain representations of the symmetric groups, and can
therefore be regarded as specializations of the symmetric functions encoding
the characters of these representations \cite{Desar}.
Classical examples include the Euler numbers \cite{Fo1}, or the Eulerian
polynomials \cite{Fo2}.
In both cases, the generating series of the relevant symmetric functions can
be obtained  as the homomorphic images of multiplicity-free series living in
some (in general noncommutative) combinatorial Hopf algebra.
For example, the generating series $\tan x + \sec x$ of the Euler numbers is
the image of the formal sum of all alternating permutations in the Hopf
algebra of free quasi-symmetric functions $\FQSym$ by its canonical character,
and its lift by Foulkes to symmetric functions is just the commutative image
of the same series \cite{NCSF1,JVNT}.

Other examples of this situation include the derangement numbers \cite{HLNT},
the numbers of parking functions \cite{NTLag}, the Abel polynomials
\cite{NTLag} or Arnold's snakes \cite{JVNT}.

These combinatorial numbers can also be dimensions of representations of the
$0$-Hecke algebra, which cannot be lifted to the generic Hecke algebra or to
the symmetric group. This is the case, for example, for linear extensions of a
poset \cite{NCSF6}, which can be directly interpreted  as free quasi-symmetric
functions.

The aim of this paper is to apply this philosophy to the Bell polynomials.
Their relation with symmetric functions is well-known, and easily extended to
noncommutative symmetric functions. At this level, it is already possible to
define a non-trivial $q$-analogue, which points toward the most promising
direction for the next step.
There is a natural choice of a representation of set partitions by
permutations in $\FQSym$ which is compatible with this $q$-analogue. The Bell
polynomials are then lifted to polynomials in noncommuting variables $Y_k$,
with coefficients in $\FQSym$. One can then consider the quasi-symmetric
functions $C_I$ which are the commutative  images of the coefficients of the
monomials $Y^I$. It turns out that they coincide with the dual immaculate
basis of \cite{BBSSZ}, up to mirror image of compositions.

The formal sum of these ``free Bell polynomials'' satisfies a simple
functional equation in terms of the dendriform structure of $\FQSym$. This
allows us to obtain expressions of the dual immaculate basis similar to (but
different from) those of Grinberg \cite{Gri}. Actually, Grinberg works
directly at the level of quasi-symmetric functions, and his formula comes in
fact from the tridendriform structure of $\WQSym$.

Finally, we discuss briefly the connection with the self-dual Hopf algebra of
set partitions introduced by Rey in \cite{rey}.

\section{Bell polynomials and noncommutative symmetric functions}

\subsection{Bell polynomials and symmetric functions}

The classical Bell (exponential) polynomials are defined by
\begin{equation}\label{eq:recB}
B_0=1, \ B_{n+1}(y_1,\ldots,y_{n+1})=\sum_{k=0}^n \binom{n}{k}B_{n-k}y_{k+1}
\end{equation}
or equivalently by the exponential generating series
\begin{equation}
\mathcal{B}(t) := \sum_{n\ge 0}B_n\frac{t^n}{n!}
                = \exp\left(\sum_{k\ge 1} y_k \frac{t^k}{k!}\right).
\end{equation}
This is reminescent of the relation between power-sums and complete symmetric
functions~\cite{Mcd}.
The power sums are
\begin{equation}
p_n = \sum_{i\ge 1} x_i^n
\end{equation}
and the complete symmetric function $h_n$ is the sum of all monomials of degree $n$, with generating series
\begin{equation}
H(t) = \sum_{n\ge 0}h_n t^n = \prod_{i\ge 1}\frac1{1-t x_i}
     = \exp\left(\sum_{k\ge 1} p_k \frac{t^k}{k}\right)
\end{equation}
so that if one sets $y_n=(n-1)!p_n$, then $B_n=n!h_n$, and for
$y_n=(-1)^{n-1}(n-1)!p_n$, $B_n$ becomes $n!e_n$, where the $e_n$ are the
elementary symmetric functions.

\subsection{Noncommutative Bell polynomials}

This is easily extended to noncommutative symmetric functions. Recall
\cite{NCSF1} that the algebra $\NCSF$ of noncommutative symmetric functions is
freely generated by noncommuting indeterminates $S_n$ playing the role of the
$h_n$. The noncommutative analogues of the $p_n$ are not canonically
determined, and one possibility is to define noncommutative power-sums of the
first kind $\Psi_n$ by orienting the Newton recursion as \cite[Prop.
3.3]{NCSF1}
\begin{equation} (n+1)S_{n+1} = \sum_{k=0}^n S_{n-k}\Psi_{k+1}.
\end{equation}
These power-sums correspond to the Dynkin elements in the descent algebras
\cite{NCSF1}.

Let $(Y_n)_{n\ge 1}$ be a sequence of noncommuting indeterminates.
For a composition $I=(i_1,\ldots,i_r)$, let $Y^I=Y_{i_1}\cdots Y_{i_r}$.
 The recurrence \eqref{eq:recB}
can be similarly oriented as
\begin{equation}
B'_{n+1}=\sum_{k=0}^n \binom{n}{k} B'_{n-k}Y_{k+1}\ 
\text{or}\ 
B''_{n+1}=\sum_{k=0}^n \binom{n}{k} Y_{k+1}B''_{n-k}.
\end{equation}
Then, $B'_n=n!S_n$ for $Y_n = (n-1)!\Psi_n$ and $B''_n=n!S_n$ for
$Y_n=(n-1)!\bar\Psi_n$ (where the $\bar\Psi_n$ are defined by the opposite
Newton recursion
\begin{equation}
(n+1)S_{n+1} = \sum_{k=0}^n \bar\Psi_{k+1}S_{n-k},
\end{equation}
and correspond to the right-sided Dynkin elements). Such noncommutative Bell
polynomials have been discussed
in \cite{SR} (see also \cite{NTshufab,ELM}).

It is enough to describe one version, since the antiautomorphism $Y_i\mapsto
Y_i$ exchanges $B'_n$ and $B''_n$.  Let us choose the second one. For now on,
we set $B=B''$.

The first values are 
\begin{equation}
\begin{split}
B_1 &= Y^1\\
B_2 &= Y^2 + Y^{11}\\
B_3 &= Y^3 + 2Y^{21} + Y^{12} + Y^{111}\\
B_4 &= Y^4 + 3Y^{31} + 3Y^{22} + 3Y^{211} + Y^{13} + 2Y^{121} + Y^{112}
       + Y^{1111} \\
B_5 &= Y^5 + 4Y^{41} + 6Y^{32} + 6Y^{311} + 4Y^{23} + 8Y^{221} + 4Y^{212}
       + 4Y^{2111} \\
    &+ Y^{14} + 3Y^{131} + 3Y^{122} + 3Y^{1211} + Y^{113} + 2Y^{1121}
     + Y^{1112} + Y^{11111} 
\end{split}
\end{equation}

The coefficients have a simple combinatorial interpretation in terms of set
partitions. We shall represent a set partition either by a set of sets or by
any sequence of its parts separated by a vertical bar.
For example, $\pi = \{\{3,4,7\},\{2,8\},\{1\},\{5,6\}\}$ will also be
represented as $347|28|1|56$.

For a set partition $\pi$ of $[n]$, let $\pi^\sharp$ be the set composition
obtained by ordering the blocks of $\pi$ w.r.t. their minima, and let $K(\pi)$
be the composition recording the sizes of the blocks of $\pi^\sharp$.
For example, if $\pi = 347|28|1|56$, then 
$\pi^\sharp=1|28|347|56$, and $K(\pi)=(1,2,3,2)$.
Then, the coefficient of $Y^I$ in $B_n$ is the number of set partitions $\pi$
of $[n]$ such that $K(\pi)=I$ (this follows by induction on $n$, as the first
block of $\pi^\sharp$ is always the one containing 1).
For example, the coefficient of $Y^{23}$ in $B_5$ enumerates the set of
partitions
\begin{equation}
12|345,\ 13|245,\ 14|235,\ 15|234,
\end{equation}
whereas the coefficient of $Y^{32}$ in $B_5$ enumerates
\begin{equation}
123|45,\ 124|35,\ 125|34,\ 134|25,\ 135|24,\ 145|23.
\end{equation}

\subsection{$q$-analogues}

In \cite[§5.1]{NCSF2}, a noncommutative analogue of the classical transformation $p_n\mapsto (1-q^n)p_n$
of ordinary symmetric functions has been defined as the algebra automorphism 
\begin{equation}\label{eq:qtrans}
\alpha_q:\ S_n(A)\mapsto S_n((1-q)A) =(1-q)\sum_{k=0}^n(-q)^kR_{1^k,n-k},
\end{equation}
(where the $R_I$ are the noncommutative ribbon Schur functions),
which is the inverse of (\cite[§ 6.1]{NCSF2})
\begin{equation}
\beta_q:\ S_n(A)\mapsto S_n\left(\frac{A}{1-q}\right)
\end{equation}
where $\frac{A}{1-q}$ denotes the (genuine) alphabet $\{q^i a_j\}$ endowed with the
order $q^ia_j<q^ka_l$ iff $i>k$ or $i=k$ and $j<l$.

In terms of the generating series of complete and elementary functions
\begin{equation}
\sigma_t(A) := \sum_{n\ge 0}t^nS_n(A)=\prod_{i\ge 1}(1-ta_i)^{-1},\quad \lambda_{-t}(A) = \sum_{n\ge 0}(-t)^n\Lambda_n(A) = \sigma_t(A)^{-1},
\end{equation}
these transformations read
\begin{equation}
\sigma_t\left(\frac{A}{1-q}\right)=\prod_{k\ge 0}^\leftarrow \sigma_{qt}(A) = \cdots\sigma_{q^2t}(A)\sigma_{qt}(A)\sigma_t(A)
\end{equation}
and
\begin{equation}
\sigma_t((1-q)A)=\lambda_{-qt}(A)\sigma_t(A).
\end{equation}

It follows from \eqref{eq:qtrans} that
the Dynkin power-sums $\Psi_n$ are limiting cases of the transformed complete
functions (see \cite{NCSF2} for an explanation of the notation)
\begin{equation}
\Psi_n = \lim_{q\rightarrow 1}\frac{S_n((1-q)A)}{1-q}
\end{equation}
and similarly
\begin{equation}
(-1)^{n-1}\Psi_n = \lim_{q\rightarrow 1}\frac{\Lambda_n((1-q)A)}{1-q}.
\end{equation}
Indeed, it is known that $\Psi_n=\sum_{k=0}^n(-1)^kR_{1^k,n-k}$, and
\begin{equation}
\sum_{n\ge 0}(-t)^n\Lambda_n((1-q)A)= (\lambda_{-qt}(A)\sigma_t(A))^{-1}=\lambda_{-t}(A)\sigma_{qt}(A)=\sigma_{qt}\left(\left(1-\frac1q\right)A\right)
\end{equation}
so that
\begin{equation}
\Lambda_n((1-q)A) = (-1)^n(q-1)\sum_{k=0}^n\left(-\frac1q\right)^nR_{1^k,n-k}.
\end{equation}

Replacing $A$ by $A/(1-q)$,
one can define $q$-analogues of $B'_n$ and $B''_n$ by
\begin{equation}
B'_n(q) = (q)_nS_n\left(\frac{A}{1-q}\right)\ \text{for $Y_n=(q)_{n-1}S_n(A)$}
\end{equation}
and
\begin{equation}
B''_n(q) = (q)_n\Lambda_n\left(\frac{A}{1-q}\right)\
\text{for $Y_n=(q)_{n-1}\Lambda_n(A)$},
\end{equation}
where $(q)_n:=\prod_{i=1}^n(1-q^i)$.

This amounts to define $B'_n(q)$ by the recursion 
\begin{equation}
B_n'(q) = \sum_{k=0}^{n-1}q^k \qbin{n-1}{k} B'_k(q)Y_{n-k}
\end{equation}
Indeed,
\begin{align}
B'_n(q)&=(q)_nS_n((\cdots+q^2A+qA)+A)\\
&= (q)_nS_n\left(\frac{qA}{1-q}+A\right)\\
&= (q)_n\sum_{k=0}^nS_k\left(\frac{qA}{1-q}\right)S_{n-k}(A)\\
&=\sum_{k=0}^{n-1}\frac{(q)_n}{(q)_{k}}q^k B'_k(q)\frac{Y_{n-k}}{(q)_{n-k-1}}+(q)_nq^nS_n \left(\frac{A}{1-q}\right)     \\
&= (1-q^n)\sum_{k=0}^n\qbin{n-1}{k}q^kB'_k(q) Y_{n-k}+q^nB'_n(q).\\
\end{align}

The same argument yields  for $B_n''(q)$
\begin{equation}
B_n''(q) = \sum_{k=0}^{n-1}q^k Y_{n-k}\qbin{n-1}{k} B'_k(q).
\end{equation}

The first values are 
\begin{equation}
\begin{split}
B'_1& = Y^1\\
B'_2& = Y^{2} + q Y^{11}\\
B'_3& = Y^{3} + q^2 Y^{21} + \left(q^2+q\right) Y^{12} + q^3 Y^{111}\\
B'_4& = Y^{4} + q^3 Y^{31} + \left(q^4+q^3+q^2\right) Y^{22} + q^5Y^{211} \\
    &   + \left(q^3+q^2+q\right) Y^{13}
      + \left(q^5+q^4\right) Y^{121} + \left(q^5+q^4+q^3\right) Y^{112}
      + q^6Y^{1111}
\end{split}
\end{equation}
and
\begin{equation}
\begin{split}
B''_1& = Y^1\\
B''_2& = Y^{2} + q Y^{11}\\
B''_3& = Y^{3} +\left(q^2+q\right) Y^{21} + q^2Y^{12} + q^3Y^{111}\\
B''_4& = Y^{4} + \left(q^3+q^2+q\right) Y^{31}
         + \left(q^4+q^3+q^2\right) Y^{22}
         + \left(q^5+q^4+q^3\right) Y^{211} \\
     & + q^3 Y^{13} + \left(q^5+q^4\right) Y^{121}
       + q^5 Y^{112} + q^6Y^{1111}.\\
\end{split}
\end{equation}
Setting $Y_i=1$ in $B'_n$ or in $B''_n$, one obtains the triangle A188919 of
\cite{Sloane}:
\begin{equation}
\begin{split}
B_1(q) &= 1\\
B_2(q) &= 1 + q\\
B_3(q) &= 1 + q + 2q^2 + q^3\\
B_4(q) &= 1 + q + 2q^2 + 4q^3 + 3q^4 + 3q^5 + q^6 \\
B_5(q) &= 1 + q + 2q^2 + 4q^3 + 7q^4 + 8q^5 + 9q^6 + 9q^7
        + 6q^8 + 4q^9 + q^{10}.\\
\end{split}
\end{equation}

The On-Line Encyclopedia of Integer Sequences \cite{Sloane} suggests that the
coefficient of $q^k$ in $B_n(q)$ should be the number of permutations of
length $n$ with $k$ inversions that avoid the dashed pattern $1-32$
\cite{Ba}.
The proof of this fact will follow from a refined interpretation for the
coefficient of $Y^I$ in $B'_n(q)$ (Proposition \ref{pr:CI} below).
Indeed,
it follows from the definition that this coefficient 
is\footnote{This is a $q$-analogue of the r.h.s. of \cite[(78)]{NTshufab}.}
\begin{equation}
\prod_{k=2}^{\ell(I)}
\left[ \genfrac{}{}{0pt}{}{i_1+\cdots+i_k-1}{i_k-1} \right]_q
  q^{i_1+\cdots+i_{k-1}}.
\end{equation}
This suggests that the counting of $1-32$-avoiding permutations
can be refined according to their descent sets.
To do this, we have to replace the binomial coefficients in the recurrence by
formal sums of permutations in an appropriate algebra.

\bigskip

As we have seen, the Bell polynomials are essentially the expressions
of the complete homogeneous symmetric functions on suitably normalized power sums.
This expression is classically given by a determinant, which has been lifted
to the noncommutative case in \cite[Theorem 2.5]{ELM} as a quasideterminant.
Expressions of the $S_n$ as quasideterminants in the $\Psi_n$ have been given
in \cite{NCSF1}, and have been extended to the $q$-analogues 
$\Theta_n(q) = (1-q)^{-1}S_n((1-q)A)$ of the $\Psi_n$ in \cite[Prop. 5.5]{NCSF2}.
A straightforward manipluation of this quasideterminant yields the following
$q$-analogue of 
\cite[Theorem 2.5]{ELM}:

\begin{proposition}\label{qdet}
$B_n'(q)$ is given by the quasideterminant
\begin{equation}
B_n'(q) = |\mathbb{B}_n'|_{1n} 
\end{equation}
where $\mathbb{B}_n'$ is the $n\times n$ matrix whose subdiagonal elements are
$-1$, the lower elements are all zero, and the elements on and above the
diagonal are
\begin{equation}
({\mathbb B}_n')_{ij}=q^{i-1}
  \left[ \genfrac{}{}{0pt}{}{n-(n-j+1)}{j-i} \right] Y_{j-i+1}\
(\text{for}\ i\le j).
\end{equation} 
\end{proposition}

\Proof
Applying the inverse transformation $\beta_q$, it is clear that
$B'_n(q)$ is also the polynomial expressing $[n]_q!S_n(A)$
in terms of the $Y_k=[k-1]_q!\Theta_k(q)$ (which are obviously a set
of free generators of $\Sym$ on $\K(q)$.
In the quasideterminant (78) of \cite{NCSF2}, expressing $S_n(A)$ in terms
of $\Theta_n(q)$ 
\begin{equation}
[n]_q \, S_n(A)
=
\left|
\begin{matrix} 
\Theta_1(q) & \Theta_2(q)    & \dots  & \Theta_{n-1}(q) & \bo{\Theta_n(q)} \\
-[1]_q      & q\,\Theta_1(q) & \dots  & q\,\Theta_{n-2}(q) & q\,\Theta_{n-1}(q)\\
 0           & - [2]_q        & \dots  & q^2\, \Theta_{n-3}(q) & q^2\, \Theta_{n-2}(q)\\
\vdots      & \vdots         & \vdots & \vdots   & \vdots \\
0           & 0              & \dots  & -[n-1]_q & q^{n-1} \, \Theta_1(q)
\end{matrix}
\right|
\ .
\end{equation}
multiply column $j$ by $[j-1]_q!$ and divide row $i$ by $[i-1]_q!$. These
operations do not change the value of the quasideterminant, except for the
last column, which yields a factor $[n-1]_q!$, so that the l.h.s. of
as \cite[Eq. (78)]{NCSF2} is now $[n]_q!S_n$. Replacing each entry
$\Theta_k(q)$ by $\frac1{[k-1]_q!}{Y_k}$ yields the desired expression.
\qed

For example,
\begin{equation}
B'_4(q) =
\left|
\begin{matrix} 
Y_1 & Y_2    & Y_3   & \bo{Y_4} \\
-1      & q\qbin{2}{0}Y_1 &  q\qbin{2}{1} Y_2  & q\qbin{2}{2}Y_3\\
 0           & -1         & q^2\qbin{3}{0}Y_1  & q^2\qbin{3}{1}Y_2 \\
0           & 0              & -1 & q^3\qbin{4}{0}Y_1
\end{matrix}
\right|
\ .
\end{equation}

Note that Equation (77) of \cite{NCSF2} allows one to express $Y_n$ in terms
of the $B_k$, as in \cite[Section 3.2.6]{ELM}.

\section{Quasi-symmetric analogues of the coefficients}

\subsection{Free quasi-symmetric functions and their dendriform structure}

Recall that for a totally ordered alphabet $A$, $\FQSym(A)$ is the algebra
spanned by the noncommutative power series~\cite{NCSF6}
\begin{equation}
\G_\sigma(A)  := \sum_{\genfrac{}{}{0pt}{}{w\in A^n}{\Std(w)=\sigma}} w
\end{equation}
where $\sigma$ is a permutation in the symmetric group $\SG_n$ and $\Std(w)$
denotes the standardization of the word $w$, {\it i.e.}, the unique
permutation having the same inversions as $w$.

All countably infinite alphabets lead to isomorphic algebras. In the sequel,
we fix one such alphabet $A$ and denote the resulting algebra $\FQSym(A)$
simply by $\FQSym$.

The multiplication rule is the same as in the Malvenuto-Reutenauer algebra
\cite{MR}, to which $\FQSym$ is therefore isomorphic\footnote{The interest of
this isomorphism resides in the fact that the Hopf algebra structure is then
induced by the ordered sum of alphabets, as for noncommutative symmetric
functions.}
\begin{equation}
\G_\alpha \G_\beta = \sum_{\gamma\in \alpha * \beta} \G_\gamma,
\end{equation}
where the convolution $\alpha*\beta$ of $\alpha\in\SG_k$ and $\beta\in\SG_l$
is defined as the sum in the group algebra of $\SG_{k+l}$ in \cite{MR}
\begin{equation}
\alpha * \beta =
\sum_{\genfrac{}{}{0pt}{}{\gamma=uv}{\Std(u)=\alpha ;\, \Std(v)=\beta}}
\gamma\,,
\end{equation}
and is interpreted here as a set. Indeed, there are no multiplicities in this
sum, and the number of terms in $\G_\alpha \G_\beta$ is $\binom{k+l}{k}$.
To reduce this number to $\binom{k+l-1}{k}$ and relate to the Bell induction
concerning $B''_n$, we can split the product according to the dendriform
structure of $\FQSym$.

The dendriform structure of $\FQSym$, originally defined in \cite{LR2}, can be
seen as inherited from the tridendriform structure of the algebra of
noncommutative formal power series over $A$ defined in \cite{NT1,NT2}.
Define bilinear operators on nonempty words $u,v$ by
\begin{eqnarray}
u\prec v=
\begin{cases}
uv &\mbox{if $\max(v) < \max(u)$}\\
0 &\mbox{otherwise},\\
\end{cases}\\
u\succ v=
\begin{cases}
uv &\mbox{if $\max(v)\geq\max(u)$}\\
0 &\mbox{otherwise}.\\
\end{cases}
\end{eqnarray}
Note that the operator $\succ$ here is $\succ+\circ$ of \cite{NT1,NT2}.
Also, our conventions differ from those of Grinberg \cite[Def. 6.1]{Gri} who
uses $\min$ instead of $\max$, and $<$ and $\ge$ instead of $>$ and $<$.
This yields
\begin{equation}
\G_\alpha \G_\beta = \G_\alpha \prec \G_\beta + \G_\alpha \succ \G_\beta\,,
\end{equation}
where
\begin{equation}\label{eq:precG}
\G_\alpha \prec \G_\beta =
\sum_{\genfrac{}{}{0pt}{}{\gamma=uv \in \alpha*\beta}%
{|u|=|\alpha| ;\, \max(v)<\max(u)}}
\G_\gamma\,,
\end{equation}
\begin{equation}
\G_\alpha \succ \G_\beta =
\sum_{\genfrac{}{}{0pt}{}{\gamma=u v\in \alpha*\beta}%
{|u|=|\alpha| ;\, \max(v)\geq\max(u)}} \G_\gamma\,.
\end{equation}
Then $x=\G_1$ generates a free dendriform algebra in $\FQSym$, isomorphic to
$\PBT$, the Loday-Ronco algebra of planar binary trees~\cite[Prop. 3.7]{LR1} and \cite[Prop. 5.3]{LR2}. 
The notations used here are those of \cite{HNT}, where more details
can be found.

In terms of the dual basis $\F_\sigma:=\G_{\sigma^{-1}}$, these operations read
\begin{align}
\F_\alpha\F_\beta &= \sum_{\gamma\in \alpha\shuffle \beta[k]}\F_\gamma,\\
\F_\alpha\prec \F_\beta &= \sum_{\gamma\in \alpha\prec \beta[k]}\F_\gamma,\\
\F_\alpha\succ \F_\beta &= \sum_{\gamma\in \alpha\succ \beta[k]}\F_\gamma,
\end{align}
where $\alpha\in\SG_k$, $\beta[k]$ denotes $\beta$ with its entries shifted by
$k$, and the $\prec$ and $\succ$ in the right-hand sides  are the half-shuffles
defined by
\begin{equation}
ua \prec vb  = (u\shuffle vb)a, \ ua\succ vb = (ua\shuffle v)b,
\end{equation}
the shuffle product $\shuffle$ being itself recursively defined by
\begin{equation}\label{eq:shuffle}
ua\shuffle vb = (ua\shuffle v)b+(u\shuffle vb)a,
\end{equation}
where $a,b$ are letters and $u,v$ words,
with the scalar 1 (representing the empty word) as neutral element for
$\shuffle$.
We shall not need dendriform products involving $1$ and thus leave these
undefined.

There is an inclusion of Hopf algebras $\Sym\hookrightarrow \PBT$ of
noncommutative symmetric functions into $\PBT$~\cite{LR1,LR2}, 
which is given
by
\begin{equation}
S_n=
(\dots((x\succ x)\succ x)\dots)\succ x\quad\text{($n$ times).}
\end{equation}
Indeed, the r.h.s is the sum of nondecreasing words, {\it i.e.},
$\G_{1\ldots n}=S_n$.

\subsection{A dendriform lift of the Bell recurrence}

We shall define the free Bell polynomials as elements of
$\K\<Y_1,Y_2,\ldots\>\otimes\FQSym$ by
\begin{equation}
\B_0=1,\quad \B_{n+1} = \sum_{k=0}^n Y_{k+1} S_{k+1} \prec \B_{n-k},
\end{equation}
the $\otimes$ sign being omitted for notational convenience, and $\prec$
is actually $\cdot\,\otimes\prec$. Recall also that $S_{k+1}=\G_{1\ldots k+1}$.
The first $\B_n$ are
\begin{equation}
\begin{split}
\B_0 & = 1 \\
\B_1 & = Y^1\G_1 \\
\B_2 & = Y^2\G_{12} + Y^{11}\G_{21} \\
\B_3 & = Y^3\G_{123} + Y^{21}(\G_{132}+\G_{231})
         + Y^{12}\G_{312} + Y^{111}\G_{321} \\
\B_4 & =
 Y^4\G_{1234} + Y^{31}(\G_{1243}+\G_{1342}+\G_{2341}) \\
 &+ Y^{22}(\G_{1423}+\G_{2413}+\G_{3412}) 
 + Y^{211}(\G_{1432}+\G_{2431}+\G_{3421}) \\ 
 &+ Y^{13}\G_{4123} + Y^{121}(\G_{4132}+\G_{4231})
 + Y^{112}\G_{4312} + Y^{1111}\G_{4321}. \\
\end{split}
\end{equation}
Clearly, for each $n\ge 0$, $\B_n$ is a (left) $\FQSym$-linear combination of
the $Y^I$ for various compositions $I$ of $n$.
The coefficient $C_I(A)$ of $Y^I$ in $\B_n$ is by definition in $\FQSym$ but
is in fact in $\PBT$:

\begin{lemma}\label{lem:CI}
Let $I=(i_1,\dots,i_r)$ be a composition. Then $C_I(A)$ is equal to $\P_T$
where $T$ is the right-comb tree whose left branches have sizes from top
to bottom equal to $i_1$, $i_2,\ldots,i_r$.
\end{lemma}

\Proof
The coefficient $C_I(A)$ is by definition equal to
\begin{equation}\label{eq:lprod}
S_{i_1}\prec (S_{i_2}\prec(\cdots (S_{i_{r-1}}\prec S_{i_r}))\cdots).
\end{equation}
Let $n$ be $i_1+\dots+i_r$ and expand that product in $\FQSym$ on the $\G$
basis. 
Consider any permutation $\tau$ such that $\G_\tau$ occurs in the expansion.
By definition of $\prec$, the first $i_1$ values are increasing
and the $i_1$-th value is $n$, then followed by a word whose standardization
$\sigma$ belongs to $S_{i_2}\prec(\cdots (S_{i_{r-1}}\prec S_{i_r}))\cdots)$.
Hence, by induction on $r$, the decreasing tree (see~\cite{LR1,HNT}) of all
these elements is obtained by grafting the shape of the decreasing tree of
$\sigma$ to the right of the root of a left branch of size $i_1$.
Conversely, since any permutation of decreasing tree $T$ belongs to the 
left product \eqref{eq:lprod}, this product is equal to $\P_T$.
\qed

\begin{example}\label{ex:Y221}{\rm
The coefficient of $Y^{221}$ in $\B_5$ is 
\begin{align*}
\G_{12}\prec (\G_{132}+\G_{231})&=
\G_{15243}+\G_{25143}+\G_{35142}+\G_{45132}\\
&+\G_{15342}+\G_{25341}+\G_{35241}+\G_{45231}\\
&=\P_{{\begin{picture}(5,5) 
\put(1,3){\circle*{0.7}}
\put(2,4){\circle*{0.7}}
\put(3,2){\circle*{0.7}}
\put(4,3){\circle*{0.7}}
\put(5,2){\circle*{0.7}}
\put(4,3){\Line(-1,-1)}
\put(4,3){\Line(1,-1)}
\put(2,4){\Line(-1,-1)}
\put(2,4){\Line(2,-1)}
\put(2,4){\circle*{1}}
\end{picture}}}
\end{align*}
}
\end{example}

For a set partition $\pi$ of $[n]$, let $\pi^\flat$ be the set composition
obtained by ordering the blocks w.r.t. their maximal values in decreasing
order, and let $K(\pi)$ be the composition recording the lengths of these
blocks. Let also $\hat\pi^\flat$ be the permutation obtained by reading the
blocks in this order, where each block is read from the smallest to the
largest element.

\begin{example}{\rm
If $\pi = 347|28|1|56$, we have
$\pi^\flat= 28|347|56|1$, so that
$\hat\pi^\flat=28347561$,
and $K(\pi)=(2,3,2,1)$.
} 
\end{example}

\begin{proposition}
\label{pr:CI}
The coefficient of $Y^I$ in $\B_n$ is the sum of all 
$\G_{\hat\pi^\flat}$
where $\pi$ ranges over set partitions such that $K(\pi)=I$:
\begin{equation}
\B_n = \sum_{\pi\vdash [n]}Y^{K(\pi)} \G_{\hat\pi^\flat}.
\end{equation}
\end{proposition}

\Proof
This follows immediately from the product rule \eqref{eq:precG}.
Indeed, as we have already seen, this coefficient is
$S_{i_1}\prec (S_{i_2}\prec(\cdots (S_{i_{r-1}}\prec S_{i_r}))\cdots)$,
which is the sum of all permutations of the form $\sigma=u_1u_2\cdots u_r$,
where each $u_j$ is an increasing word of length $i_j$ whose last letter is
greater than all letters to its right.
\qed

These permutations are not those avoiding the pattern $1-32$ as above, but
those avoiding the pattern $21-3$. As observed in \cite[Prop. 1]{Cla}, there
is a statistic-preserving bijection between both classes (in this case, the
Sch\"utzenberger involution, which preserves the inversion number). 
Since the inverse major index and the inversion number are equidistributed on
the set of permutations having a decreasing tree of a given shape
\cite[Theorem 2.7]{HNT2}, the generating polynomial of $\inv$  on $1-32$ or of $\imaj$ on
$21-3$-avoiding permutations with descent composition $I$ is
\begin{equation}
c_I(q)=C_I\left(\frac1{1-q}\right)
\end{equation}
where the alphabet $\frac1{1-q}$ (the principal specialization of $\FQSym$) is defined by
\begin{equation}
A\mapsto \{q^n|n\ge 0\} \ \text{ordered by $q^i<q^j$ iff $i>j$}.
\end{equation} 
Then,
\begin{equation}
\G_\sigma\left(\frac1{1-q}\right)=F_I\left(\frac1{1-q}\right)
\end{equation}
where $I$ is the descent composition of $\sigma^{-1}$, and
the specializations of the fundamental quasi-symmetric functions are
\begin{equation}
F_I\left(\frac1{1-q}\right)
=\frac{q^{\maj(I)}}{(q)_n}.
\end{equation}
where $\maj(I)=\maj(\tau)$ for any permutation $\tau$ of shape $I$.
Thus, specializing the $\FQSym$ coefficients in $\B_n$, we obtain 
\begin{equation}
(q)_n\B_n\left(\frac1{1-q}\right) = B''_n(q).
\end{equation}

\bigskip

Now that we understand that $c_I(q)$ is the principal specialization of a
quasi-symmetric function, we can replace it by the commutative image $C_I(X)$
of $C_I(A)$ in $\QSym$. Recall that if the letters $a_i$ of our underlying alphabet $A$ 
 are replaced by commuting
variables $x_i$, then $\F_\sigma$ becomes the fundamental quasi-symmetric
function $F_I(X)$ where $I$ is the descent composition of $\sigma$.

\begin{example}
{\rm Recording the recoil compositions of the permutations occuring in
$C_{221}$, we find
\begin{equation}
C_{221}(X)=
F_{1121} + F_{1211} + F_{122} + F_{131} + F_{212} + 2F_{221} + F_{311}.
\end{equation}
This may be compared with the dual immaculate basis of \cite{BBSSZ}:
\begin{equation}
\SG_{221}^*=
F_{1121} + F_{113} + F_{1211} + 2F_{122} + F_{131} + F_{212} + F_{221}.
\end{equation}
}
\end{example}

\begin{theorem}\label{th:dualimm}
Define the bar involution on $QSym$ by $\overline{F_I} = F_{\bar I}$,
where for $I=(i_1,\ldots,i_r)$, $\bar I =(i_r,\ldots,i_1)$ denotes the mirror
composition. Then, the dual immaculate basis
is given by
\begin{equation}
\SG_I^* = \overline{C_I(X)}.
\end{equation}
\end{theorem}

\Proof
According to \cite[Prop. 3.37]{BBSSZ},
\begin{equation}
\SG_I^* = \sum_T F_{D(T)} 
\end{equation}
where the sum runs over all standard immaculate tableaux of shape $I$, and
$D(T)$ denotes the descent composition of $T$ as defined in \cite{BBSSZ}.

This should not be confused with the usual descent composition of a permutation,
which will be denoted by $C(\sigma)$.

A standard immaculate tableau $T$ of shape $I$ is a planar representation
of a set partition $\pi=(\pi_1,\ldots,\pi_r)$, whose blocks have been ordered in such a way
that $\min(\pi_1)<\min(\pi_2)<\ldots<\min(_pi_r)$, and such that $|\pi_j|=i_j$.

For example, the standard immaculate tableaux of shape $I=(2,2,1)$ are
\footnotesize
\begin{equation}
\young{1 & 5\cr 2& 4\cr 3\cr}\ \ \
\young{1 & 5\cr 2&3 \cr4 \cr}\ \ \
\young{1 &4 \cr 2&5 \cr 3\cr}\ \ \
\young{1 & 4\cr 2& 3\cr 5\cr}\ \ \
\young{1 &3 \cr2 &5 \cr 4 \cr}\ \ \
\young{1 & 3\cr 2& 4\cr 5\cr}\ \ \
\young{1&2\cr3&5\cr 4\cr}\ \ \
\young{1 & 2\cr 3& 4\cr 5\cr}
\end{equation}
\normalsize
The descent composition $D(T)$ encodes the set $\{i|i+1\ \text{is in a lower row}\}$,
which is therefore the recoil set of the permutation $\hat T$ obtained by reading the
rows of $T$ from bottom to top and from left to right.

The descent compostions $D(T)$ of the above tableaux are
\begin{equation}
113,\ 122,\ 1121,\ 131,\ 122,\ 1211,\ 212,\ 221.
\end{equation}
These are the recoil compositions  of the permutations $\hat T$ 
\begin{equation}
32415,\ 42315,\ 32514,\ 52314,\ 42513,\ 52413,\ 43512,\ 53412,
\end{equation}
whose descent compositions are obviously always $\bar I=(1,2,2)$.

The Sch\"utzenberger involution $\nu$ sends a permutation  $\sigma\in\SG_n$
to the permutation $\sigma'=\nu(\sigma)$ obtained by replacing each entry $i$
of $\sigma$ by $n+1-i$ and then reading the resulting word from right to left.
In other words, $\sigma'= \omega\sigma\omega$, where $\omega=n\,n-1\cdots 21$.

The effect of $\nu$ on descent compositions is 
the mirror image $C(\sigma')=\overline{C(\sigma)}$. Thus, for any standard
immmaculate tableau $T$ of shape $I$ and descent composition $J$, applying
$\nu$ to the permutation $\sigma=\nu(\hat T)$ yields
\begin{equation}
C(\sigma) = C(\nu(\hat T)) =\overline{C(\hat T)} = \overline{\bar I} = I
\end{equation}
while
\begin{equation}
C(\sigma^{-1})=\overline{D(\hat T)} = \bar J.
\end{equation}
Moreover, $\nu$ exchanges the maxima and the minima, so that the rows of the ribbon diagram of $\sigma$
form the blocks of a set partition $\pi=(\pi_1,\ldots,\pi_r)$ ordered in such a way that 
$\max(\pi_1)>\max(\pi_2)>\ldots>\max(\pi_r)$, and $|\pi_j|=i_j$. That is, $\sigma=\hat\pi^\flat$ for a
set partition such that $K(\pi)=I$.
Therefore,
\begin{equation}
\SG_I^* = 
\sum_{K(\pi)=I}\overline{F_{K(\hat\pi^\flat)}}=
\overline{C_I(X)}.
\end{equation}

On our running example, the permutations $\nu(\hat T)$ are
\begin{equation}
15243,\ 15342,\ 25143,\ 25341,\ 35142,\ 35241,\ 45132,\ 45231.
\end{equation}
Their recoil compositions are, in order
\begin{equation}
311,\ 221,\ 1211,\ 131,\ 221,\ 1121,\ 212,\ 122.
\end{equation}
\qed

\begin{note}{\rm
Theorem \ref{th:dualimm} provides an expression of the dual immaculate basis
very similar to the one obtained by Grinberg \cite{Gri}. However,
our dendriform operations are different. The relations between both
constructions will be clarified in the forthcoming section.
}
\end{note}

\begin{corollary}
The $C_I(X)$ form a basis of $QSym$, so that the $\P_T$ indexed by right combs
form a section of the projection $\PBT\rightarrow QSym$.
\end{corollary}

\begin{corollary}
$c_I(q)$ is given by the Bj\"orner-Wachs $q$-hook-length formula \cite{BW}:
the inversion polynomial of the set of permutations having
a decreasing tree of shape $T$  is given by the
same hook length formula as for the inverse major index, which is
\begin{equation}
\sum_{{\mathcal T}(\sigma)=T}q^{l(\sigma)}
=
[n]_q! \prod_{v\in T}\frac{q^{\delta_v}}{[h_v]_q}\,,
\end{equation}
where $v$ runs over the vertices of $T$, $h_v$ is the number of vertices
of the subtree with root $v$, and 
 $\delta(v)$ is the number of vertices in the right subtree of $v$.
\end{corollary}

See \cite{HNT2} for a short proof of the version used here. 

\begin{example}{\rm
For $I=221$, the hook-lengths $h_v$ of the right comb are $5,3,1,1,1$, and the
cardinalities of the right subtrees are $3,1$. Hence,
\begin{equation} 
	c_{221}(q) = [5]_q!\frac{q^{3+1}}{[5]_q[3]_q[1]_q^3}
={q}^{4}+2\,{q}^{5}+2\,{q}^{6}+2\,{q}^{7}+{q}^{8}.
\end{equation}
}
\end{example}

Applying the Sch\"utzenberger involution $\nu$ to $C_I$  amounts to
sending $\G_\sigma$ to $\check\G_\sigma= \G_{\omega\sigma\omega}$, so that we can also reformulate Theorem \ref{th:dualimm} as
\begin{equation}
\check C_I(X)=\SG_I^*.
\end{equation}
\begin{corollary}[\cite{BBSSZ2}]
$\overline{\SG_I^*}$ is the characteristic of an indecomposable $0$-Hecke
algebra module.
\end{corollary}

\Proof
The inverses of the permutations occuring in $C_I(A)$ are the linear
extensions of a poset (a binary tree in this case), hence form the basis of a
$0$-Hecke module, see Section 3.9 of~\cite{NCSF6}. Thus, the commutative
image $C_I(X)=\overline{\SG_I^*}$ of $C_I(A)$ in $QSym$ is the characteristic
of this module. The same is true of their images by the Sch\"utzenberger
involution, the poset being now a binary tree turned upside-down.
Moreover, these modules are indecomposable, since they are of the form
${\bf N}_\sigma$ described in \cite[Definition 4.3]{NCSF6}. Indeed, the right-comb
tree associated with a composition $I$ is the shape of the devrasing tree of the maximal permutation
$\omega(I)$ of the descent class $I$. This permutation (which is self-inverse)
spans the (one-dimensional) socle of the module, which must be contained in
any submodule, so that no submodule can be a direct summand.
\qed

See \cite{NCSF4} or \cite{Thib} for an introduction to the representation theory of the 0-Hecke algebra.

\section{The dual immaculate basis}

\subsection{Half-shuffles and $\FQSym$}

The shuffle product on $\K\<A\>$ can be recursively defined by
\eqref{eq:shuffle}:
\begin{equation}
ua \shuffle vb  = (ua\shuffle v)b + (u\shuffle vb)a
\end{equation}
or symmetrically by
\begin{equation}\label{eq:halfsh2}
au \shuffle bv  = a(u\shuffle bv) + b(au\shuffle v),
\end{equation}
where $a,b\in A$ and $u,v\in A^*$.
The half-shuffles are also known as chronological products. Both ways
of splitting the shuflle
can be used to define a dendriform structure on $\FQSym$.
In this section, we shall use the second possibility \eqref{eq:halfsh2} and set
\begin{equation}
au\prec' bv = a(u\shuffle bv),\ \text{and}\ au\succ' bv = b(au\shuffle v).
\end{equation}
If $\gamma$ is a linear combination of some permutations $\rho$,
we write for short $\F_\gamma$ for the same linear combination of the
$\F_\rho$.
On $\FQSym$, we set, for $\sigma\in\SG_k$ and $\tau\in\SG_l$
\begin{equation}
\F_\sigma\prec' \F_\tau = \F_{\sigma\prec' \tau[k]},\ \text{and}\
\F_\sigma\succ' \F_\tau = \F_{\sigma\succ' \tau[k]}.
\end{equation}
Here, $\succ'$ coincides with Grinberg's $\succ$, which we will denote by $\succ_G$ to avoid
confusions.
\begin{lemma}
For any two words $u$ and $v$,
\begin{equation}
\label{eq:halfsh}
u\prec' v =  \sum_{v_1v_2=v} (-1)^{|v_1|} (\overline{v_1}u) \shuffle v_2,
\end{equation}
where $\bar w$ denotes the mirror image of a word $w$.
\end{lemma}
\Proof By definition
\begin{equation}
au\prec' bv = au\shuffle bv - au\succ' bv,
\end{equation} 
and since 
\begin{equation}
au\succ' bv = bau\prec' v,
\end{equation}
the result follows by induction.
\qed

For example,
\begin{equation}
1234\prec' 567 =
1(234 \shuffle 567) = 1234 \shuffle 567 - 51234 \shuffle 67 + 651234 \shuffle 7 - 7651234.
\end{equation}
In $\FQSym$, this implies
\begin{equation}
\label{eq:hsfqs}
\F_\sigma\prec' \F_\tau =
\sum_{uv = \tau[k]}(-1)^{|u|}\F_{(\bar u \cdot\sigma)\shuffle v}. 
\end{equation}
For example,
\begin{equation}\label{eq:exdend}
\F_{2143}\prec'\F_{312}
=
\F_{2143\shuffle 756}-\F_{72143\shuffle 56}+\F_{572143\shuffle 6}-\F_{6572143}.
\end{equation}

\subsection{Descents in half-shuffles}

For $w\in A^*$, let $\alph(w)\subseteq A$ be the set of letters occuring in
$w$.
The following property appears as Lemma 4.1 in \cite{NTsuper}:
\begin{lemma}
\label{desshuf}
If $\alph(u)\cap\alph(v)=\emptyset$, then 
\begin{equation}
\<u\shuffle v\>=\<u\>\<v\>\,,
\end{equation}
where the linear map $\< \,\>$ is defined by $\<u\>:=F_{C(u)}$, where $C(u)$
denotes the descent composition of the word $u$.
In particular, the descents of the elements of a shuffle on disjoint alphabets
depend only on the descents of the initial elements.
\end{lemma}

There is a refined  statement for the dendriform
half-products~\cite[Theorem 4.2]{NTsuper} that we adapt to our half-products:

\begin{theorem}
\label{deshalf}
Let $u=u_1\cdots u_k$ and $v=v_1\cdots v_\ell$ be two nonempty words of respective lengths $k$ and
$\ell$.
If $\alph(u)\cap\alph(v)=\emptyset$, then
\begin{equation}
\<u\prec' v\>=\<\sigma\prec' \tau\>
\end{equation}
where $\sigma=\std(u)$ and $\tau=\std(v)[k]$ if $u_1<v_1$,
and $\sigma=\std(u)[\ell]$ and $\tau=\std(v)$ if $u_1>v_1$.
\end{theorem}

For example, we have
\begin{equation}
\begin{split}
\<14\prec' 23\> &= \<1423+1243+1234\> = F_{22} + F_{31} + F_4 \\
\<12\prec' 34\> &= \<1234+1324+1342\> = F_4 + F_{22} + F_{31},
\end{split}
\end{equation}
whereas
\begin{equation}
\begin{split}
\<24\prec' 13\> &= \<2413+2143+2134\> = F_{22} + F_{121} + F_{13} \\
\<34\prec' 12\> &= \<3412+3142+3124\> = F_{22} + F_{121} + F_{13}.
\end{split}
\end{equation}

\subsection{Projection onto $QSym$}

Let $\pi:\ \FQSym\rightarrow QSym$ be the canonical projection sending
$\F_\sigma$ to $F_{C(\sigma)}$. Then $\pi$ is compatible with the products of
both structures, and 
we can define half-products $\prec'$ and $\succ'$ on $QSym$ in such a way that
$\pi$ will be compatible with the half-products as well.
This follows from an even
more general trivial property. Let us consider $\sigma$ and $\sigma'$ having
the same descents and $\tau$ and $\tau'$ also. Any word $w$ in the shifted shuffle
$\sigma\ssh\tau=\sigma\shuffle \tau[k]$ ($\sigma\in\SG_k$) is characterized by the sequence of positions of the letters
of $\sigma$. Then the word $w'$ in 
$\sigma'\ssh\tau'$ with the same sequence
of positions of the letters of $\sigma'$ has obviously same descents as $w$.
Since it is true for all elements, it is in particular true when one sums over
subsets of the shuffle, \emph{e.g.}, $\prec'$ and $\succ'$.
Thus, if $C(\sigma)=I$ and $C(\tau)=J$, we can define in $QSym$
\begin{equation}\label{eq:QSdend}
F_I\prec' F_J
:=\pi(\F_\sigma\prec' \F_\tau),
\end{equation}

For example,
\begin{equation}
\F_{21}\prec' \F_{132} = \F_{21354}+\F_{23154}+\F_{23514}+\F_{23541},
\end{equation}
so that
\begin{equation}
\label{eq:exgauF}
F_{11}\prec' F_{21} = F_{131}+ F_{221}+ F_{32}+ F_{311},
\end{equation}
which could have been computed as well from
\begin{equation}
\F_{21}\prec' \F_{231} = \F_{21453}+\F_{25153}+\F_{24513}+\F_{24531}.
\end{equation}

Applying Lemma~\ref{desshuf}, we can project
\eqref{eq:hsfqs} to $QSym$.
If the descent composition of $u$ is $H$ and that of $\sigma$ is $I$, the
descent composition of $\bar u$ is the conjugate composition $H^\sim$, and if
$u_1>\sigma_1$, the descent composition of $\bar u \cdot \sigma$ is
$H^\sim\cdot I$.
We have therefore
\begin{equation}
F_I\prec' F_J =
\sum_{J\in\{HK,H\triangleright K\}}(-1)^{|H|}F_{ H^{\sim}\cdot I}F_K,
\end{equation}
where for two compositions $H=(h_1,\ldots,h_r)$ and $K=(k_1,\ldots,k_s)$
\begin{equation}
HK = (h_1,\ldots,h_r,k_1,\ldots,k_s)\ \text{and}\ H\triangleright K=(h_1,\ldots,h_r+k_1,\ldots,k_s).
\end{equation}

For example,
\begin{equation}
\F_{21}\prec' \F_{132}
= \F_{21\shuffle 354} - \F_{321\shuffle 54} + \F_{5321\shuffle 4}
- \F_{45321},
\end{equation}
so that
\begin{equation}
\pi(\F_{21}\prec' \F_{132})
= F_{11} F_{21} - F_{111} F_{11} + F_{1111} F_{1} - F_{2111},
\end{equation}
corresponding to the decompositions
$J=\emptyset\cdot 21$, $1\triangleright 11$, $2\cdot 1$,
and $21 \cdot\emptyset$.

Since the antipode $S$ of $QSym$ is given by 
\begin{equation}
S(F_H)=(-1)^{|H|}F_{ H^{\sim}},
\end{equation} 
and the coproduct by 
\begin{equation}
\Delta(F_J)=\sum_{J\in\{HK,H\triangleright K\}}F_H\otimes F_K
\end{equation}
we obtain in this way \cite[Detailed version (ancillary file), Theorem 3.15]{Gri}:
\begin{theorem}
\label{th:precs}
Let the $\prec'$ product on $QSym$ be defined by \eqref{eq:QSdend}. Then,
for $f,g\in QSym$,
\begin{equation}
\label{eq:fprecg}
f\prec' g = \sum_{(g)}\left(S(g_{(1)})\bullet f\right)g_{(2)}
\end{equation}
where $F_I\bullet F_J := F_{I\cdot J}$ and
$\Delta g = \sum_{(g)} g_{(1)}\otimes g_{(2)}$.
\end{theorem}

For example, applying $\pi$ to \eqref{eq:exdend}, we obtain
\begin{equation}\label{eq:exFF}
F_{121}\prec' F_{12}
=
F_{121}F_{12}-F_{1121}F_{2}+F_{2121}F_1-F_{12121.}
\end{equation}

\subsection{Quasi-differential operators}

Theorem~\ref{th:precs} can be reformulated in terms of quasi-differential
operators, as in~\cite{Gri}.
According to \cite[Definition 4.5]{NCSF2}, for $f\in QSym$,
\begin{equation}
f(X-Y) = \sum_I S(F_I)(Y)\cdot R_I^\perp(f)(X),
\end{equation}
where the quasi-differential operator $R_I^\perp$ is defined as in
\cite{Gri}, \emph{i.e.}, as the adjoint of the linear map $f\mapsto R_I f$.
Indeed, recall that $f(X-Y)$ is defined as $f((-Y)\hat+ X)$, where $\hat+$
is the ordinal sum of alphabets,
$f(-Y)=S(f)(Y)$.
Now,
\begin{align}
f(X\hat+ Y) &= \sum_{J,K}\<\Delta f, R_J\otimes R_K\>F_J(X)F_K(Y)\\
&= \sum_{J,K}\<f, R_J R_K\>F_J(X)F_K(Y)\\
&= \sum_{J,K}\<R_J^\perp f R_K\>F_J(X)F_K(Y)\\
&= \sum_J F_J(X)\sum_K\<R_J^\perp f R_k\>)F_K(Y)\\
&= \sum_J F_J(X)(R_J^\perp f)(Y).
\end{align}
Replacing $X,Y$ by $-Y,X$, this becomes
\begin{equation}
f(-Y\hat+ X)= \sum_J F_J(-Y)(R_J^\perp f)(X)= \sum_J S(F_J)(Y)\cdot R_J^\perp(f)(X),
\end{equation}

We can therefore rewrite \eqref{eq:fprecg} as
\begin{equation}
f\prec' g (X) = [g(X-Y)\bullet_Y f(Y)]_{Y=X}.
\end{equation}
For example
\begin{align}
F_{12}(X-Y)
&=F_{12}((-Y)\hat+ X) \\
&=F_{12}(-Y)+F_{11}(-Y)F_1(X)+F_1(-Y)F_2(X)+F_{12}(X)\nonumber\\
&=-F_{12}(Y)+F_2(Y)F_1(X)-F_1(Y)F_{2}(X)+F_{12}(X).
\end{align}
Taking the $\bullet$ product with $F_{121}(Y)$, we obtain
\begin{equation}
F_{12}(X-Y)\bullet F_{121}(Y)
= -F_{12121}(Y)+F_{2121}(Y)F_1(X)-F_{1121}(Y)F_{2}(X)+F_{121}(Y)F_{12}(X),
\end{equation}
and setting $Y=X$, we recover \eqref{eq:exFF}.

Similarly, we have for the right product
\begin{equation}
u\succ' v = \sum_{u_1u_2=u}(-1)^{|u_1|}u_2\shuffle \overline{u_1}v.
\end{equation}
Thus, on $QSym$,
if one defines $\btr$ by $F_I\btr F_J=F_{I \triangleright J}$, 
then
\begin{equation}
f\succ'g = \sum_{(f)}\left(S(f_{(1)})\btr g\right) f_{(2)}.
\end{equation}
This is precisely \cite[Theorem 3.7]{Gri}.
For example,
\begin{equation}
\F_{132}\succ' \F_{21}
= \F_{132\shuffle 54 - 32\shuffle 154 + 2\shuffle 3154 - 23154},
\end{equation}
which projects onto
\begin{equation} 
F_{21}\succ' F_{11} = F_{21}F_{11} - F_{11}F_{21} + F_{1}F_{121} - F_{221}
   = F_{131} + F_{1121} + F_{122} + F_{1211}. 
\end{equation}

\subsection{Standard dendriform structures}

With the usual definitions
\begin{equation}
ua\prec vb = (u\shuffle vb)a,\quad ua\succ vb = (ua\shuffle v)b,
\end{equation}
the half-shuffle identity becomes
\begin{equation}
u\prec v = \sum_{v_1v_2=v}(-1)^{|v_2|}u\overline{v_2}\shuffle v_1
\end{equation}
which induces an operation on $QSym$ defined by
\begin{equation}
F_I\prec F_J = \sum_{J\in\{HK,H\triangleright K\}}
                 (-1)^{|K|}F_{I\triangleright K^\sim}F_H,
\end{equation}
or equivalently,
\begin{equation}
f\prec g = \sum_{(g)}g_{(1)}\left(f\btr S(g_{(2)})\right).
\end{equation}

For example,
\begin{equation}
\F_{21}\prec \F_{132}=\F_{35421}+\F_{35241}+ \F_{32541}+\F_{23541}
\end{equation}
and
\begin{equation}
F_{11}\prec F_{21}=F_{2111}+F_{221}+F_{1211}+F_{311}.
\end{equation}
Then,
\begin{equation}
F_I\prec F_J = F_I(Y)\btr_Y F_J(X-Y)|_{Y=X}.
\end{equation}
For example,
\begin{equation}
F_{21}(X-Y) = F_{21}(X) - F_{1}(Y)F_{2}(X) + F_{2}(Y)F_{1}(X) - F_{21}(Y)   
\end{equation}
and
\begin{equation}
F_{11}(Y)\btr_Y F_{21}(X-Y) = F_{11}(Y)F_{21}(X) - F_{12}(Y)F_{2}(X)
                           +  F_{13}(Y)F_{1}(X) -  F_{131}(Y)  
\end{equation}
so that
\begin{equation}
F_{11}\prec F_{21}= F_{11}F_{21} - F_{12}F_{2}+ F_{13}F_{1} - F_{131}.
\end{equation}

For the right product, we have
\begin{equation}
u\succ v = \sum_{u_1u_2=u}(-1)^{|u_2|}u_1\shuffle v\overline{u_2},
\end{equation}
and on $QSym$,
\begin{equation}
f\succ g = \sum_{(f)}f_{(1)}\left(g\bullet S(f_{(2)})\right).
\end{equation}
For example,
\begin{equation}
\F_{21}\succ \F_{231} = \F_{21\shuffle 453 - 2\shuffle 4531 + 45312}
\end{equation}
which yields on $QSym$
\begin{equation}
F_{11}\succ F_{21}=F_{11}F_{21}- F_{1} F_{211}+ F_{212} 
=
F_{131}+F_{122}+ F_{221}+ F_{32}+ F_{1121}+ F_{212}. 
\end{equation}

\subsection{Grinberg's operations}

In \cite{Gri}, a left product $\prec_G$ on $QSym$ is induced from an operation
on monomials.
This operation is the commutative image of the left tridendriform product on
words defined in \cite{NTpark} (with the use of $\min$ instead of $\max$, so
that we shall consistently denote it by $\prec'$), and amounts to taking the
canonical projection $\pi: \WQSym\rightarrow QSym$ of the tridendriform
product of $\WQSym$:
\begin{equation}
M_I\prec_G M_J := \pi(\M_u\prec' \M_v)
\end{equation}
where $u,v$ are any packed words of evaluation $I,J$.

For example, to evaluate $F_{11}\prec_G F_{21}$, we compute
$M_{11}\prec_G (M_{21}+M_{111})$, hence
$\M_{21}\prec'(\M_{132}+\M_{121})$:
\begin{equation*}
\M_{21}\prec'\M_{132}
=
\M_{21243}+
\M_{21354}+
\M_{31243}+
\M_{41243}+
\M_{51243}+
\M_{41253}+
\M_{31254},
\end{equation*}
\begin{equation*}
\M_{21}\prec' \M_{121}=
\M_{21232}+
\M_{31212}+
\M_{21343}+
\M_{31242}+
\M_{41232}
\end{equation*}
so that
\begin{equation}
\begin{split}
F_{11}\prec_G F_{21} &= 3M_{1211}+2M_{1121}+M_{1112}+M_{122}+M_{131}+4M_{1111}\\
                     &= F_{131}+F_{122}+F_{1121}+F_{1211},
\end{split}
\end{equation}
which is different from all previous examples such as \eqref{eq:exgauF}.

This operation can then be lifted to $\FQSym$ as
\begin{equation}
\F_\sigma \prec_G \F_\tau := \F_{\sigma[l]\prec'\tau},
\end{equation}
where $l=|\tau|$, and $au\prec' bv$ is defined as $a(u\shuffle bv)$.
Then Equation~\eqref{eq:halfsh} yields \cite[Theorem 3.7]{Gri}

\begin{equation}
F_I\prec_G F_J =
 \sum_{J\in\{HK,H\triangleright K\}}(-1)^{|H|}F_{ H^{\sim}\triangleright I}F_K.
\end{equation}

For example,
\begin{equation}
\F_{21}\prec_G \F_{132} = \F_{54132}+ \F_{51432}+ \F_{51342}+ \F_{51234},
\end{equation}
and one can check that
\begin{equation}
54\prec' 132 = 54\shuffle 132-154\shuffle 32+3154\shuffle 2 -23154,
\end{equation}
so that
\begin{equation}
F_{11}\prec_G F_{21} = F_{11}F_{21}-F_{21}F_{11}+F_{121}F_1-F_{221}.
\end{equation}

Alternatively, we can describe $\prec_G$ as
\begin{equation}
\F_\sigma\prec_G \F_\tau = \F_\tau\succ' \F_\sigma.
\end{equation}

On $QSym$, this translates as
\begin{equation}
f\prec_G g = \sum_{(g)} g_{(2)}\left( S(g_{(1)})\btr f\right)
\end{equation}
which is now \cite[Theorem 3.7]{Gri} in its original form.

Grinberg's expression of the dual immaculate basis can now be restated as
\begin{equation}
\SG_I^* = (\cdots (F_{i_r}\succ' F_{i_{r-1}})\succ'\cdots)\succ' F_{i_1}.
\end{equation}

For example,
\begin{equation}
\SG_{221}^*=(F_1\succ' F_2)\succ' F_2= (F_{12}+F_{21})\succ' F_2.
\end{equation}

\begin{note}
On $\FQSym$, the operation $\prec_G$ is not a left
dendriform product in the usual sense, as it is in fact a flipped right
product.
The $\prec'$ operation on $\WQSym$ which induces it does not preserve the
standard $\FQSym$ subalgebra of $\WQSym$.
\end{note}

Finally, it is also possible to deduce Theorem \ref{th:dualimm} from Grinberg's results.
The following argument has been suggested by an anonymous referee.

Let $f\mapsto \bar f$ be the antinvolution of $\FQSym$ defined by $\overline{\F_\sigma}=\F_{\omega\sigma\omega}$,
and denote by $\pi:\ \FQSym\rightarrow QSym$ the canonical projection $\pi(\F_\sigma)=F_{C(\sigma)}$.
As already observed, $\pi(\bar f)=\overline{\pi(f)}$. It is easy to see that
\begin{equation}\label{eq:abba1}
\overline{a\prec b}=\bar b\succ_G\bar a\quad\text{for all $a,b\in\FQSym$}.
\end{equation}
Moreover, in \cite[§6]{Gri}, it is said that
\begin{equation}\label{eq:abba2}
\pi(a)\prec_G \pi(b) = \pi(a\succ b) \quad\text{for all $a,b\in\FQSym$}.
\end{equation}
Now, let $I=(i_1,i_2,\ldots,i_r)$ be a composition. Then, \cite[Corollary 4.7]{Gri} shows that
\begin{align*}
\SG_I^* &= h_{i_1}\prec_G (  h_{i_2}\prec_G(\cdots (   h_{i_r}\prec_G 1)\cdots))\\
&= h_{i_1}\prec_G (  h_{i_2}\prec_G(\cdots (h_{i_{r-1}}\prec_G   h_{i_r})\cdots))\\
&= \pi(S_{i_1}) \prec_G(  \pi(S_{i_2}) \prec_G(\cdots\prec_G   \pi(S_{i_{r-1}}) \prec_G \pi(S_{i_r}))\cdots))\\
& \text{(since $h_k=\pi(S_k)$ for all $k$)}\\
&=\pi(((\cdots (S_{i_r} \succ_G S_{i_{r-1}}) \succ_G\cdots) \succ_G S_{i_2}) \succ_G S_{i_1} )\\
& \text{(by applying \eqref{eq:abba2} many times)}\\ 
&=\pi(((\cdots (\overline{S_{i_r}} \succ_G \overline{S_{i_{r-1}}}) \succ_G\cdots) \succ_G \overline{S_{i_2}}) \succ_G \overline{S_{i_1}} )\\
&\text{(since $\overline{S_k}=S_k$)}\\
&=\pi\left(\overline{S_{i_1}\prec(S_{i_2}\prec(\cdots\prec(S_{i_{r-1}}\prec S_{i_r})\cdots))}\right)\\
&\text{(by applying \eqref{eq:abba1} many times)}\\
&=\pi\left(\overline{C_I(A)}\right) \quad\text{as seen in the proof of Lemma \ref{lem:CI}}\\
&=\overline{\pi(C_I(A))}\\
&= \overline{C_I(X)}. 
\end{align*}

Hence, Theorem \ref{th:dualimm} follows from \cite[Corollary 4.7]{Gri}.

\section{Hopf algebras of set partitions}

We have seen that the Bell polynomial $\B_n$ can be identified with the formal
sum of permutations avoiding $21-3$ (up to the $Y^I$ which can be
reconstructed from the descent sets). This raises the question of the
existence of a Hopf subalgebra or quotient of $\FQSym$ whose bases are
naturally labeled by these permutations.

The most obvious Hopf algebra of set partitions is $\WSym$, or symmetric
functions in noncommuting variables (not to be confused with noncommutative
symmetric functions $\NCSF$).
In \cite{HNT1}, a quotient of $\WSym$ isomorphic to $\NCSF$  and a $QSym$
subalgebra of its dual are related to Bell polynomials. In \cite{BCLM},
analogues of the Bell polynomials in various other Hopf algebras are
considered.

\subsection{The Bell Hopf algebra}

The Hopf algebra $\WSym$ is cocommutative. Quite often, combinatorial objects
also admit a self-dual Hopf algebra structure.
Such an algebra has been constructed by Rey \cite{rey,rey2} for set
partitions from the Burstein-Lankham correspondence, a combinatorial
construction derived from the patience sorting algorithm. More precisely,
the Bell classes of Rey are indexed by permutations avoiding $23-1$, and these
are their minimal elements (for the weak order). We can modify the
construction so as to have classes whose maximal elements avoid $21-3$ as
follows.

Let $A$ be a totally ordered alphabet. The (modified) Bell congruence on $A^*$
is generated by the relations%
\footnote{As observed by Grinberg (private communication), an alternative
description is $bcu\equiv buc$ for $b<c$ if all letters of $u$ are $<b$.}
\begin{equation}
buca\equiv buac\quad\text{if $a<b<c$ and all letters of $u$ are smaller than $b$.}
\end{equation}

This is a refinement of the reverse sylvester congruence, which is defined by
the same relations without restriction on $u$.

As for the sylvester or reverse sylvester congruences, we have:

\begin{theorem}
\label{th:posets}
The Bell equivalence classes of permutations,
({\it i.e;}, equivalence classes of words with no equal letters under the Bell congruence)
 are intervals of the right weak
order on the symmetric group. These intervals consist of the linear extensions
of posets which will be explicitly described below, and the maximal elements
of these intervals are the $21-3$-avoiding permutations.
\end{theorem}

Identifying set partitions $\pi$ with their representatives $\hat\pi^\flat$ as $21-3$-avoiding
permutations, we set
\begin{equation}
P_\pi=\sum_{\sigma\equiv\pi} \F_\sigma.
\end{equation}

\begin{theorem}[Rey \cite{rey}]\label{th:rey}
The $P_\pi$ span a Hopf subalgebra of $\FQSym$.
\end{theorem}

\subsection{The poset of set partitions}

Since the long version of Rey's paper \cite{rey} has never been published, we
shall provide detailed proofs of his results (in our modified version), to
which we add the explicit description of the posets.

\subsubsection{From permutations to set partitions}
Let us first define an insertion algorithm. This is a rewriting of the
Patience Sorting Algorithm defined by Burstein and Lankham adapted to our
setting, which amounts to applying some trivial involutions on words.
Precisely, if we invert the total order of the alphabet, our algorithm becomes
Algorithm 3.1 of \cite{rey2}, our blocks being Rey's piles read downwards.

The first object that we create is a set partition, that we shall regard as
ordered in such a way that the maximal elements of the blocks decrease from
left to right.

Let $w=w_1\dots w_n$ be a word with no repeated letters over a totally ordered
alphabet.
Put $S=\emptyset$. Then, read $w$ from left to right and for each letter $w_i$
do
\begin{itemize}
\item Step 1: Let $s$ be the block of $S$ whose maximal element is smaller
than (or equal to) $w_i$ and is greater than all other  maximal elements
smaller than $w_i$.
\item Step 2: If $s$ does not exist, 
add $\{w_i\}$ to $S$. 
If $s$ exists, insert $w_i$ into it.
\end{itemize}
The result of this algorithm will be denoted by $PSA(w)$.

We shall display the blocks of $S$ as columns, increasing from top to bottom
(so that an element is inserted at the bottom of its column), the columns
being ordered from left to right with their maximal elements in decreasing
order.

For example, the insertion of $(3,1,2,6,4,5,7)$ follows the steps:

\medskip

$S =\emptyset$
$\to$
$S = \smalltableau{3}$
$\to$
$S = \smalltableau{3 & 1}$
$\to$
$S = \smalltableau{3 & 1\\\ & 2}$
$\to$
$S = \smalltableau{3 & 1\\ 6 & 2}$
$\to$
$S = \smalltableau{3 & 1\\ 6 & 2 \\\ & 4}$
$\to$
$S = \smalltableau{3 & 1\\ 6 & 2 \\\ & 4\\\ & 5}$
$\to$
$S = \smalltableau{3 & 1\\ 6 & 2 \\ 7 & 4\\\ & 5}$

\subsubsection{From set partitions to posets}
Now, starting from a set partition $S$, we shall build a poset $P$ on the ground set $\{1,2,\ldots n\}$
whose relation $>_P$ is generated by the following requirements: order
$S$ as above, and, for any element $x$ of a block $s$ of $S$, write $x>_Px'$
where $x'$ is the greatest element smaller than $x$ in $s$ and write $x>_Px''$
where $x''$ is the smallest element greater than $x$ in the 
block immediately to the left of $s$ in $S$ (if such a block and such an element exist).
The transitive closure of $>_P$ defines a poset which  will be denoted by $P(S)$ or by $P(w)$ if $S$ is the
result of the (modified) patience sorting algorithm applied to $w$.
Note that by removing the edges of $P(S)$ where the element above is greater
than the element below, one recovers $S$ itself. Hence it makes sense to
refer to the columns of $S$  as the corresponding columns
of the poset $P(S)$. For $w$ in $P$, we shall write $C(w)$ for its column.

As usual with posets, we shall only represent the covering relations.
For example, starting with the set $S$ computed before, we get (orienting
the Hasse diagrams upside-down, with minimal elements at the top):

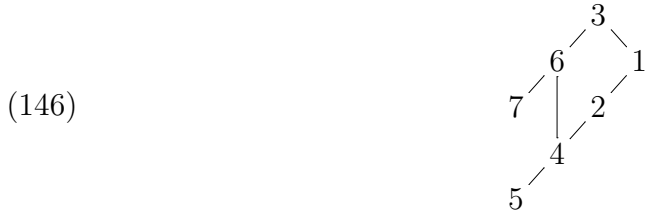
\entrymodifiers={+<4pt>}
\begin{equation}
\vcenter{\xymatrix@C=2mm@R=2mm{
*{}   &  *{} & {3}\ar@{-}[dr]\ar@{-}[dl] \\
*{}   & {6}\ar@{-}[dd]\ar@{-}[dl]  & *{} & {1}\ar@{-}[dl] \\
{7}   & *{}  & {2}\ar@{-}[dl] \\
*{}   & {4}\ar@{-}[dl] \\
{5} \\
      }}
\end{equation}

Indeed, the first column of $S$ is $(3,6,7)$ and $6>3$ and $7>6$.
Now, on the second column we have $5>4>2>1$ and the relations between these
elements and the first column are: $1>_P3$, $2>_P3$, $4>_P6$, and $5>_P6$, whence the
Hasse diagram above.

\medskip
Starting with $S=(6,10,11|2,4,8,9|3,7|5|1)$, one gets
$S = \smalltableau{6&2&3&5&1\\ 10&4&7 \\ 11&8 \\\ & 9}$
and

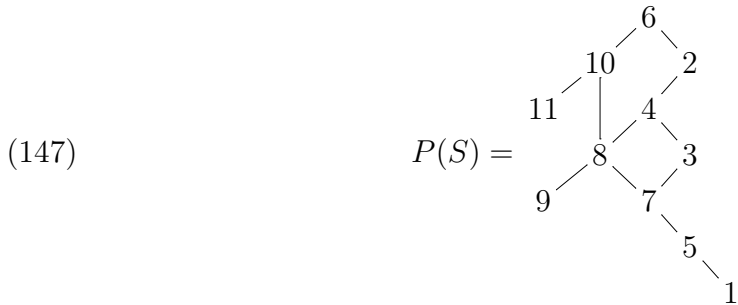
\entrymodifiers={+<4pt>}
\begin{equation}P(S)=
\vcenter{\xymatrix@C=2mm@R=2mm{
*{}   &  *{} & {6}\ar@{-}[dr]\ar@{-}[dl] \\
*{}   & {10}\ar@{-}[dd]\ar@{-}[dl]  & *{} & {2}\ar@{-}[dl] \\
{11}   & *{}  & {4}\ar@{-}[dl]\ar@{-}[dr] \\
*{}   & {8}\ar@{-}[dl]\ar@{-}[dr] & *{} & {3}\ar@{-}[dl] \\
{9} & *{} & 7\ar@{-}[dr] \\
*{} & *{} & *{} & 5\ar@{-}[dr] \\
*{} & *{} & *{} & *{} & 1\\
      }}
\end{equation}

Note that in this representation, the columns of $P(S)$ are the straight lines
going from south-west to north-east.

\subsection{Proofs}

We shall now prove Theorem \ref{th:posets}:
all elements in a given Bell class have the same poset $P(w)=P(PSA(w))$, and
the linear extensions (taken from top to bottom in our representation) of such
a poset give elements of the same Bell class. Hence being in the same Bell
class is equivalent to having the same poset.

\medskip
\subsubsection{All elements in a given Bell class have the same poset}

First, let us prove that all elements in a given Bell class have the same
poset, or, equivalently, the same result by the PSA algorithm. We only need to
prove this for two elements $w$ and $w'$ obtained from one another by a single
rewriting rule. So, let us consider $w=vbuac$ and $w'=vbuca$ where all letters
of $u$ are smaller than $b$. When $b$ is inserted, it is the bottommost element of
its column, and it remains so during the insertion of all $u$, since all the
letters of $u$ are smaller than $b$. Now, a letter $c>b$ is necessarily inserted into
a column to the left of $C(b)$ or into it, independently of the insertion of
$a$.
Meanwhile, a letter $a$ smaller than $b$ is necessarily inserted into a column
strictly to the right of $C(b)$ (whether or not $c$ has already been
inserted). Thus, the insertions of $c$ and $a$ do not interfere with each
other.
So $w$ and $w'$ satisfy $PSA(w)=PSA(w')$.

Note that the condition that all letters of $u$ are smaller than $b$ is
necessary to ensure that $c$ and $a$ cannot both be inserted in the same
column, which would prevent $w$ and $w'$ from yielding the same result by the
modified patience sorting algorithm.

\medskip
\subsubsection{All linear extensions of a poset are Bell congruent}

Let us now prove that all linear extensions of a poset $P(\sigma)$ are Bell
congruent to $\sigma$.%
\footnote{Grinberg (private communication) has proposed another proof, an
induction of the length of $\sigma$.}
First note that the poset obtained from a permutation $\sigma$ has a unique
minimal element $m$ at its top: indeed, it is easy to see that whenever
$p<q$, each element of the $q$th column is $>_P$ to the smallest element of
the $p$th column. 
Now, $\sigma$ is a linear extension of it, since when we insert
a new letter $r$, we do not get a new covering relation of the form $s>_Pr$.

One can thus rebuild the corresponding set partition $PSA(\sigma)$ since 
each column of $P(\sigma)$ is a saturated chain of the poset $P(\sigma)$,
and conversely, if we remove all edges of the form $u>_P v$ satisfying
$u<v$ from the Hasse diagram of $P(\sigma)$, then only these saturated chains remain.
Let us define a {\it Bell poset} as a poset $P(S)$ for a set partition $S$
of some finite set of integers, not necessarily $\{1,2,\ldots,n\}$.
It is easy to see that removing a maximal element from a Bell poset yields a Bell
poset, which implies that any down-closed poset of a Bell poset is again a Bell poset.
So, let us now consider two values $a<c$, both maximal elements of $P$, and let
$P'$ be the poset obtained from $P$ by removing $a$ and $c$.
By induction on $n$, we can assume that all linear extensions of $P'$ are
Bell-congruent with one another. To prove that the same holds for $P$, we only
need to prove that there exists one linear extension $v$ of $P'$ such that
$v\cdot ac\equiv v\cdot ca$: indeed, all linear extensions of a poset can be reached
from one another by such exchanges.

First, note that $C(a)$ is to the right of $C(c)$ since they are both maximal
elements of their column. Let us prove that there exists a linear extension of
$P'$ that ends with a letter $b$ such that $a<b<c$, possibly followed by
letters smaller than $a$.
We shall use twice a very simple observation: if $x$ is the last letter of a
column, then there are linear extensions of the associated poset where $x$ is
only followed by letters smaller than itself. Indeed, the letters that must
appear after $x$ in all linear extensions of the poset are letters belonging
to columns to the right of $C(x)$, hence all smaller than $x$.

Now, if $C(c)$ is not the column immediately to the left of $C(a)$, define
$b$ as the maximal element of this column. It satisfies $a<b<c$ and all
elements after it in the linear extensions of $P'$ are smaller than $b$, and
in fact smaller than $a$. So there are linear extensions of $P'$ where the
rewriting $wca\equiv wac$ is possible.
Otherwise, since $a$ and $c$ are both maximal elements of $P$, they cannot be
comparable.
Define $b$ as the letter immediately above $c$ in its column. It exists since
otherwise $a$ would be connected to $c$ (definition of the posets from the
sets of columns). Now, $b$ might not be a maximal element of $P'$ but in any
case, the values that are below $b$ in $P'$ are strictly to its right
in the partition, hence are all smaller than $a$, again. So there are linear
extensions of $P'$ ending with $b$ and then letters smaller than $a$.
So in both cases, we found a value $b$ allowing to rewrite one linear
extension $w$ of $P'$ concatenated to $ac$ with $w \cdot ca$, thus proving
that all linear extensions of a poset $P(\sigma)$ are Bell congruent to
$\sigma$.

\medskip
\subsubsection{Our posets are regular}

Now that we have obtained the equivalence between Bell classes and the linear
extensions of our posets, it only remains to prove that the posets
(regarded as labeled posets)  are regular
to conclude that the classes are intervals of the weak order \cite[Theorem 5.8]{BW}.
We shall make use of a very simple property: if a value $u$ belongs to a
column to the left of a value $v$, and if $u<v$ then $u<_P v$. Indeed,
consider the column $C'$ immediately to the left of $C(v)$.
Then $v$ is $>_P$ the smallest value $v'$ greater than itself
in $C'$, hence below all values smaller than itself. Now the same property
holds for $v'>v$ on the next column to its left and so on.

Let us now consider two elements $x$ and $z$ such that $x<_P z$, so that
$C(x)$ is equal to, or to the left of $C(z)$, which
we shall write as $C(x)\preceq C(z)$.

. Two cases may appear: $x<z$ or $x>z$.
Let us consider the first case and let $y$ be another element of $P$ such that
$x<y<z$. In that case, either $C(y)\preceq C(z)$ or $C(y)\succeq C(x)$.
 In both cases, apply our previous result and conclude
that either $x>_Py$ or $y<_Pz$.
Let us now consider the second case where $x>z$ and let $y$ be another element
of $P$ such that $x>y>z$. If $C(y)\preceq C(z)$ or $C(y)\succeq C(x)$, 
the same argument as before applies. Otherwise, it means
that  $C(x)\preceq C(y)\preceq C(z)$. 
But in that case, since there is a
saturated chain from $x$ down to $z$, all elements in between are
either comparable to $x$ (if they are below the path) or to $z$ (if they are
above the path). So in all cases, we proved that $y$ satisfies either $x<_P y$
or $y<_Pz$, hence proving that our posets are regular.

\medskip
\subsubsection{The maximal elements of the Bell classes avoid $21-3$}

Let us finally prove that the maximal elements of the Bell classes avoid the
pattern $21-3$. There are several ways to prove this. One could say for
example that the Bell classes are in bijection with set partitions, hence
splitting the permutations of size $n$ into as many classes as the number of
permutations avoiding $21-3$. Now, if a permutation $\sigma$ does not avoid
$21-3$, let us consider such a pattern with letters $ba\dots c$ in $\sigma$
where the distance between $a$ and $c$ is minimal. Then the letter immediately
to the left of $c$ (if it exists: the other case is easy)
has to be smaller than $b$, hence denoted naturally by
$a'$.  All letters $x$ between $b$ and $a'$ are necessarily smaller than $b$,
otherwise
$ba\cdots x\cdots a'c$ would have $ba\cdot x$ as a shorter pattern $21-3$, which would
contradict the minimality of the initial one. So we are in the conditions of
the Bell congruence, and  $\sigma$ is not a maximal element of its class. So
all maximal elements necessarily avoid the pattern $21-3$. Since both sets,
the maximal elements of Bell classes and the permutations avoiding $21-3$ have
same cardinality, they must be equal.

Another proof consists in applying the greedy algorithm that takes a Bell
poset as entry and takes the linear extension where at each step one chooses
the maximal available value. It is easy to see that these elements necessarily
avoid the pattern $21-3$ and conclude in the same way as before.

Note that the minimal elements are not characterized by pattern avoidance.

\medskip
\subsubsection{Proof of Theorem~\ref{th:rey}}

As for Theorem \ref{th:rey}, it follows from a well-known result
about algebras defined by congruences, see \emph{e.g.}, Theorem 2.1
of~\cite{NT16} or Chapter~2.1 of~\cite{JNZ}.

\subsection*{Acknowledgements} This research has been partially supported by
the ANR program CARMA.
Thanks also to Darij Grinberg for his useful comments on the first
version of this text.

\end{document}